\documentclass{amsart}
\usepackage{amsmath, amsthm, amssymb, amscd, epsfig}
\usepackage{pinlabel}

\begin{document}

\newtheorem{theorem}{Theorem}[section]
\newtheorem{lemma}[theorem]{Lemma}
\newtheorem{proposition}[theorem]{Proposition}
\newtheorem{corollary}[theorem]{Corollary}
\newtheorem{conjecture}[theorem]{Conjecture}
\newtheorem{question}[theorem]{Question}
\newtheorem{problem}[theorem]{Problem}
\newtheorem*{claim}{Claim}
\newtheorem*{criterion}{Criterion}
\newtheorem*{stability_thm}{Stability Theorem~\ref{stability_theorem}}

\theoremstyle{definition}
\newtheorem{definition}[theorem]{Definition}
\newtheorem{construction}[theorem]{Construction}
\newtheorem{notation}[theorem]{Notation}

\theoremstyle{remark}
\newtheorem{remark}[theorem]{Remark}
\newtheorem{example}[theorem]{Example}

\numberwithin{equation}{subsection}

\def\Z{\mathbb Z}
\def\H{\mathbb H}
\def\RP{\mathbb{RP}}

\def\L{\mathcal L}

\def\length{\textnormal{length}}
\def\arg{\textnormal{arg}}
\def\Li{\textnormal{Li}}

\def\Id{\textnormal{Id}}
\def\PSL{\textnormal{PSL}}
\def\til{\widetilde}

\title{Bridgeman's orthospectrum identity}
\author{Danny Calegari}
\address{Department of Mathematics \\ Caltech \\
Pasadena CA, 91125}
\date{\today}

\begin{abstract}
We give a short derivation of an identity of Bridgeman concerning orthospectra of 
hyperbolic surfaces.
\end{abstract}

\maketitle

\section{Introduction}

In \cite{Bridgeman}, Martin Bridgeman proves a beautiful identity concerning orthospectra of 
hyperbolic surfaces with totally geodesic boundary. Let $\Sigma$ be a hyperbolic surface 
with totally geodesic boundary. An {\em orthogeodesic} 
is a geodesic segment properly immersed in $\Sigma$, which is perpendicular to 
$\partial \Sigma$ at its endpoints. The set of orthogeodesics is countable, and their 
lengths are proper. Denote these lengths by $l_i$ (with multiplicity). Bridgeman's identity is:

\begin{theorem}[Bridgeman, \cite{Bridgeman}]
With notation as above,
$$\sum_i \L(1/\cosh^2(l_i/2)) = -\pi^2\chi(\Sigma)/2$$
\end{theorem}

where $\L$ is the Rogers' dilogarithm function (to be defined below). 

Treating the function $\L$ as a black box for the moment, 
the identity has the form 
$$\sum_i \ell(l_i) = \text{a term depending only 
on the topology of } \Sigma$$ 
The proof is very, very short and elegant. By the Gauss--Bonnet 
theorem, the term on the right is equal to $1/8$ of the volume of the unit tangent bundle of 
$\Sigma$. Almost every tangent vector on $\Sigma$ can be exponentiated to a geodesic on $\Sigma$ 
which intersects the boundary in finite forward and backward time (by ergodicity of the 
geodesic flow on a closed hyperbolic surface obtained by doubling; see e.g. \cite{Katok_Hasselblatt}). 
If $v$ is such a tangent vector, and $\gamma_v$ is the associated geodesic arc, then $\gamma_v$ is
homotopic keeping endpoints on $\partial \Sigma$ to a unique orthogeodesic (which is the unique 
length minimizer in its relative homotopy class). 

\begin{figure}[htpb]
\labellist
\small\hair 2pt
\pinlabel $\omega$ at 50 57
\endlabellist
\centering
\includegraphics[scale=1]{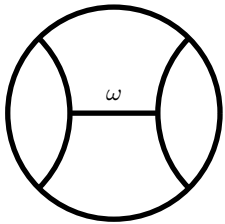}
\caption{An orthogeodesic $\omega$ lifts to the crossbar of a letter H}\label{circles_figure_3}
\end{figure}

The volume of the set of $v$ associated to 
a given orthogeodesic $\omega$ can be computed as follows. Lift $\omega$ to the universal cover, 
together with lifts of the boundary geodesics it ends on. The three together make a letter ``H''
in which $\omega$ is the crossbar; see Figure~\ref{circles_figure_3}. 
Any $\gamma_v$ in the homotopy class of $\omega$ lifts to a unique geodesic segment in the
universal cover with endpoints on the sides of the H. Therefore the volume of the set of such $v$
depends only on the geometry of the H, which in turn depends only on $l=\length(\omega)$. This
volume is $8\ell(l)$ with notation as above.

To complete the proof of Bridgeman's identity therefore, it suffices to show $\ell(l)=\L(1/\cosh^2(l/2))$.
Bridgeman derives this in several pages of calculations. The purpose of this note is
to give a short derivation of this fact, using elementary hyperbolic geometry.

\section{Derivation of $L$}

The ``ordinary'' polylogarithms $\Li_k$ can be defined by a Taylor series
$$\Li_k(z) = \sum_{n=1}^\infty \frac {z^n} {n^k}$$
which converges for $|z|<1$, and extends by analytic continuation. 
Taking derivatives, one sees that they satisfy the identities
$$\Li_k'(z) = \Li_{k-1}(z)/z$$ 
thereby giving rising to recursive integral formulae for these functions. The special case 
$\Li_0(z)$ is the familiar geometric series for $z/(1-z)$, so $\Li_1(z) = -\log(1-z)$ and
$$\Li_2(z) = -\int \frac {\log(1-z)} {z} dz$$
The Rogers dilogarithm is a ``normalization'' of $\Li_2$ given by the 
formula 
$$\L(z) = \Li_2(z) + \frac 1 2 \log(|z|)\log(1-z)$$ for real $z<1$. 
One sees that the Rogers dilogarithm is obtained by symmetrizing the integrand 
for the integral expression for $\Li_2$ under the involution $z \to 1-z$:
$$\L'(z) = -\frac {1}{2} \left(\frac {\log(1-z)}{z} + \frac {\log(z)}{1-z} \right)$$
A classic reference for this material is Lewin \cite{Lewin}.

\medskip

\begin{figure}[htpb]
\labellist
\small\hair 2pt
\pinlabel $L$ at 23 50
\pinlabel $R$ at 77 50
\pinlabel $U$ at 50 77
\pinlabel $D$ at 50 23
\pinlabel $\infty$ at 85 15
\pinlabel $0$ at 85 85
\pinlabel $x$ at 15 85
\pinlabel $1$ at 15 15
\endlabellist
\centering
\includegraphics[scale=1]{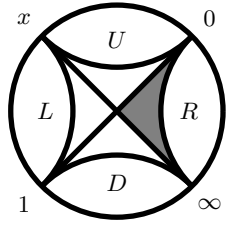}
\caption{An H spans an ideal quadrilateral, which is dissected into four semi-ideal triangles}\label{circles_figure_2}
\end{figure}

We now explain how to compute the volume of the set of vectors $v$ tangent to a
geodesic $\gamma_v$ intersecting the left and right sides of the H associated to an orthogeodesic
$\omega$. The four ideal vertices of the H span an ideal quadrilateral $Q$. The diagonals
of this quadrilateral subdivide it into four semi-ideal triangles. We denote the left and
right sides of the H as $L$ and $R$, and the other two edges of $Q$ as $U$ and $D$. Similarly,
denote the four triangles $T_L$, $T_R$, $T_U$ and $T_D$ labeled according to which edge of
the quadrilateral they bound; see Figure~\ref{circles_figure_2} (the triangle $T_R$ is colored
gray in the figure). The crossbar of the H (corresponding to the orthogeodesic $\omega$ itself) 
is not depicted in the figure.

We identify the ideal circle with $\RP^1$ in such a way that the vertices of $Q$
are (in circular order) $ 0,x,1,\infty$, where $ \infty,0$ are the ideal vertices of the gray triangle. 
Call $\alpha$ the  (hyperbolic) angle of the gray triangle at its vertex. By elementary hyperbolic
trigonometry, $x = (1+\cos(\alpha))/2 = \tanh^2(l/2)$ where $l$ is the distance between 
$L$ and $R$ (i.e. the length of the given orthogeodesic). We use the parameters $l$, $x$ and $\alpha$
interchangeably in the sequel.
We will compute $\ell$ implicitly as a function of $x$, and show that it is a multiple 
of the Rogers dilogarithm function, thus verifying Bridgeman's identity.

\medskip

Every vector $v$ in $Q$ exponentiates to a (bi-infinite) geodesic $\gamma_v$, 
and we want to compute the volume of the set of vectors $v$ for which the 
corresponding geodesic intersects both $L$ and $R$. The point of the decomposition in the
figure is that for $v$ in $T_L$ (say), the geodesic $\gamma_v$ intersects $L$ 
whenever it intersects $R$, so we only need to compute the volume of the $v$ in $T_L$ 
for which $\gamma_v$ intersects $R$. Similarly, we only need to compute the volume of the $v$ 
in $T_R$ for which $\gamma_v$ intersects $L$. For $v$ in $T_U$, we compute the volume 
of the $v$ for which $\gamma_v$ does {\em not} intersect $U$ (since these are exactly the 
ones that intersect both $L$ and $R$), and similarly for $T_D$.

The crux of the matter is that these volumes can be expressed in terms of 
integrals of simple harmonic functions. Let $\chi_L$ denote the harmonic function 
on the disk which is $1$ on the arc of the circle bounded by $L$, and $0$ on the 
rest of the circle. This function at each point is equal to $1/2\pi$ times 
the visual angle (i.e. the length in the unit tangent circle) subtended by 
the given arc of the circle, as seen from the given point in the hyperbolic plane. 
Define $\chi_R$, $\chi_U$ and $\chi_D$ similarly. Then the total volume we need 
to compute is equal to
$$4\pi \left( \int_{T_L} 2\chi_R + \int_{T_U} (1 - 2\chi_U) \right)$$
(here we have identified $\int_{T_L} \chi_R = \int_{T_R} \chi_L$ by symmetry, 
and similarly for the other pair of terms). Let us approach this a bit more systematically. 
We introduce three functions $A(\cdot)$, $B(\cdot)$ and $C(\cdot)$ as follows.
If $\alpha$ denotes as above the angle at the nonideal vertex of triangle $T_R$, we define 
$$\int_{T_R} \chi_R = A(\alpha), \quad \int_{T_R} \chi_U = B(\alpha), \quad \text{and} \quad
\int_{T_R} \chi_L = C(\alpha)$$

The integral we want to evaluate can be expressed easily in terms of explicit 
rational multiples of $\pi$, and the functions $A,B,C$. These functions satisfy obvious identities:
$$C(\alpha) = \int_{T_R} 1 - A(\alpha) - 2B(\alpha) = \pi-\alpha - A(\alpha) - 2B(\alpha)$$
and
$$A(\alpha) + B(\pi - \alpha) = \pi/3$$
where the last identity comes by observing that we are integrating a certain function 
over an ideal triangle, and observing that the average of this function under the 
symmetries of the ideal triangle is equal to the constant function $1/3$. 
In particular, we see that we can express everything in terms of $A$. 
After some elementary reorganization, we see that the contribution $V(\alpha)$ 
to the volume of the unit tangent bundle of the surface associated to this 
particular orthogeodesic is
$$V(\alpha) = \pi^2(8 - 16/3) - 4\pi\alpha - 8\pi(A(\alpha) - A(\pi - \alpha))$$

\begin{figure}[htpb]
\labellist
\small\hair 2pt
\pinlabel $0$ at 51 3
\pinlabel $\alpha$ at 82 43
\pinlabel $x$ at 77 3
\endlabellist
\centering
\includegraphics[scale=1]{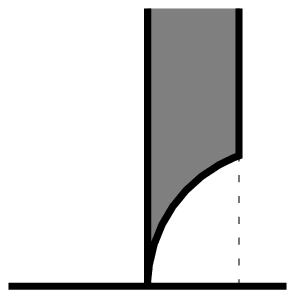}
\caption{}\label{circles_figure}
\end{figure}

It remains to actually compute $A(\alpha)$. To do this it makes sense to move to the 
upper half-space model, and move the endpoints of the interval to $0$ and $\infty$. 
The harmonic function is equal to $1$ on the negative real axis, and $0$ on the 
positive real axis. It takes the value $\theta/\pi$ on the line $\arg(z) = \theta$. 
The area form in the hyperbolic metric is proportional to the Euclidean area form, 
with constant $1/\text{Im}(z)^2$. In other words, we want to integrate 
$\arg(z)/\pi\, \text{Im}(z)^2$ over the region indicated in figure~\ref{circles_figure}, 
where the nonideal angle is $\alpha$, and the base point is $0$.

If we normalize so that the circular arc is part of the semicircle from $0$ to $1$, 
then the real projection of the vertical lines in the figure are $0$ and $x$. 

There is no elementary way to evaluate this integral, so instead we evaluate its 
derivative as a function of $x$ where as before, $x = (1+\cos(\alpha))/2$. 
This is the definite integral
$$A'(x) = \int_{y = \sqrt{x-x^2}}^\infty \frac{\tan^{-1}(y/x)}{\pi y^2} dy$$
Integrating by parts givesÊ
$$A'(x) = \frac {\alpha}{\pi\sin{\alpha}} + 
\frac 1 \pi \int_{y = \sqrt{x-x^2}}^\infty \frac x {y(y^2+x^2)} dy$$
This evaluates to
$A'(x) = (\alpha/\pi\sin{\alpha}) - 1/\pi ( \log(1-x)/2x)$.
Thinking of $V(\alpha)$ as a function of $x$, we get
$$V'(x) = -4\pi d\alpha/dx - 8\pi(A'(x) + A'(1-x)) = -8\L'(x)$$
Comparing values at $ x=0$ we see that $ V(x)=8\L(1-x)$ (where we use the
symmetry of $\L'$ under $x \to 1-x$). Substitute $x=\tanh^2(l/2)$ 
and the identity is proved.

\begin{remark}
The paper \cite{Dupont_Sah}
by Dupont and Sah relates Rogers dilogarithm to volumes of $\text{SL}(2,\mathbb{R})$-simplices, 
and discusses some interesting connections to conformal field theory and lattice model 
calculations. They cite an older paper of Dupont for the explicit calculations; these are somewhat 
tedious and unenlightening; however, he does manage to show that the Rogers dilogarithm is 
characterized by the Abel identity. In other words, Dupont shows:
\begin{lemma}[Dupont \cite{Dupont}, Lemma~A.1]
Let $f:(0,1) \to \mathbb{R}$ be a three times differentiable function satisfying
$$f(s_1) - f(s_2) + f\left(\frac{s_2}{s_1}\right) - f\left(\frac{1-s_1^{-1}}{1-s_2^{-1}}\right) 
+ f\left(\frac{1-s_1}{1-s_2}\right)=0$$
for all $0 < s_2 < s_1 < 1$. 
Then there is a real constant $\kappa$ such that $f(x) = \kappa \L(x)$ where $\L(x)$ is the 
Rogers dilogarithm (up to an additive constant).
\end{lemma}

Our geometric argument can be reformulated in homological algebraic terms, though since its
virtue is its simplicity, we have not pursued this.
\end{remark}

\section{Acknowledgments}
Danny Calegari was supported by NSF grant DMS 0707130. I would like to thank the
anonymous referee for useful comments, and Martin Bridgeman
for a beautiful talk at Caltech in which he explained his identity. I would also like
to thank Hidetoshi Masai for catching some errors in an earlier version of this paper.

\end{document}